\documentclass[leqno]{article}

\usepackage[T1]{fontenc}
\usepackage[utf8]{inputenc}

\usepackage{amsthm}
\usepackage{amssymb}
\usepackage{amsmath}

\usepackage{graphicx}
\usepackage{float}
\usepackage{caption}
\usepackage{url}
\usepackage{hyperref}

\newtheorem{theorem}{Theorem}[section]
\newtheorem{prop}[theorem]{Proposition}
\newtheorem{lemma}[theorem]{Lemma}
\newtheorem{coro}[theorem]{Corollary}

\author{M.~Hellus\footnote{Fakult\"at f\"ur Mathematik, Universit\"at Regensburg, Germany, michael.hellus@mathematik.uni-regensburg.de}
\title{On the Frobenius number of certain numerical semigroups}\and A.~Rechenauer\footnote{Anton Rechenauer, antonrechenauer@gmail.com}
\title{On the Frobenius number of certain numerical semigroups}\and R.~Waldi\footnote{Fakult\"at f\"ur Mathematik, Universit\"at Regensburg, Germany, rolf.waldi@mathematik.uni-regensburg.de}}
\title{On the Frobenius number of certain numerical semigroups}

\begin{document}

\maketitle

\begin{abstract}Let $0<\lambda\leq1$, $\lambda\notin\left\{\frac24, \frac27, \frac2{10}, \frac2{13}, \ldots\right\}$, be a real and $p$ a prime number, with $[p,p+\lambda p]$ containing at least two primes. Denote by $f_\lambda(p)$ the largest integer which cannot be written as a sum of primes from $[p,p+\lambda p]$. Then
\[f_\lambda(p)\sim\left\lfloor2+\frac2\lambda\right\rfloor\cdot p\text{, as }p\text{ goes to infinity.}\]
Further a question of Wilf about the 'Money-Changing Problem' has a positive answer for all semigroups of multiplicity $p$ containing the primes from $[p,2p]$. In particular, this holds for the semigroup generated by all primes not less than $p$. The latter special case was already shown in a previous paper.
\end{abstract}

\vspace{.3cm}

\noindent \textit{MSC 2020}: 11D07; 11P32; 20M14.

\noindent \textit{Keywords}: Numerical semigroups, Diophantine Frobenius problem, Wilf's conjecture on numerical semigroups, Goldbach conjecture.

\section{Introduction}

A \textit{numerical semigroup} is an additively closed subset $S$ of $\mathbb N$ with $0\in S$ and only finitely many positive integers outside from $S$, the so-called \textit{gaps} of $S$. The \textit{genus} $g$ of $S$ is the number of its gaps. The set $E=S^*\setminus(S^*+S^*)$, where
$S^*=S\setminus\{0\}$, is the (unique) minimal system of generators of $S$. Its elements are called the \textit{atoms} of $S$; their number $e$ is the \textit{embedding dimension} of $S$. The \textit{multiplicity} of $S$ is the smallest element $p$ of $S^*$.

From now on we assume that $S\neq\mathbb N$. Then the greatest gap $f$ is called the \textit{Frobenius number} of $S$.

For a certain class of numerical semigroups we shall study the relationship between the various invariants mentioned above.

In particular for some of these semigroups we will give an affirmative answer to

\vspace{.1cm}

\noindent\textbf{Wilf's question} \cite{wilf_1978}: Is it true that
\[\tag{1}\frac g{1+f}\leq\frac{e-1}e\text{ ?}\]
We shall consider the following semigroups: Let $p$ be a prime, $\lambda$ a positive real number, $I_\lambda(p)$ the interval $[p,p+\lambda p]$ and $D_\lambda$ the set of all primes $p$ such that $I_\lambda(p)$ contains at least two primes. For such a $p$ we denote by $S_\lambda(p)$ the numerical semigroup generated by all primes from $I_\lambda(p)$ and by $f_\lambda(p)$ its Frobenius number. According to Bertrand's postulate, $D_1$ is the set $\mathbb P$ of all primes, further $\mathbb P\setminus D_\lambda$ is finite for all $\lambda>0$ by the prime number theorem.

Let $p_1=2$, $p_2=3$, $p_3=5$, \ldots be the sequence of prime numbers in natural order and let $S_n$ be the semigroup generated by all primes not less than $p_n$. Proposition \ref{prop11}\,b) below generalizes the corresponding assertation \cite[Proposition~5]{hellus_rechenauer_waldi_19} about $S_n$.

\vspace{.1cm}

\begin{prop}

\label{prop11}Let $\lambda>0$, $p\in D_\lambda$ and $S$ any numerical semigroup of multiplicity $p$ containing $S_\lambda(p)$.

\begin{enumerate}

\item[a)] There is an integer $C(\lambda)>0$ such that for $p>C(\lambda)$, the semigroups $S$ from above satisfy Wilf's inequality (1).

\item[b)] In case $\lambda=1$ formula (1) holds for all p.
In particular (1) is true for $S=S_n$.\hfill$\qedsymbol$

\end{enumerate}

\end{prop}

This will be seen in section~\ref{section_Wilf_certain_num_sem_grps}.

In section \ref{section_upper_bounds}, we shall show the following result.

\vspace{.1cm}

\begin{theorem}

\label{theorem12}Let $0<\lambda\leq1$, $\lambda\notin\{\frac24,\frac27,\frac2{10},\frac2{13},\ldots\}$, be a real. Then 
\[\lim_{p\in D_\lambda\atop p\to\infty}\frac{f_\lambda(p)}p=\left\lfloor2+\frac2\lambda\right\rfloor.\]
For $\lambda=\frac2m$ with an integer $m\geq2$ and $m\equiv1\mod3$, at least
\[\pushQED{\qed}\limsup_{p\to\infty}\frac{f_\lambda(p)}p=\left\lfloor2+\frac2\lambda\right\rfloor.\qedhere\popQED\]

\end{theorem}

In particular, since $f_\lambda(p)$ is decreasing as a function of $\lambda$, for each $\lambda>0$ there is a constant $c_\lambda$, such that $f_\lambda(p)\leq c_\lambda\cdot p$ for all primes $p\in D_\lambda$, cf. \cite[Remark~2\,c)]{hellus_rechenauer_waldi_19}.

According to our table \texttt{t4\_quotient\_not\_always\_6.pdf} from \cite{repo} possibly the series $\frac{f_{\frac12}(p)}p$ may not converge as $p$ goes to infinity. Hence we do not expect that $\frac{f_{\frac2m}(p)}p\sim 2+m$ holds in the exceptional cases $m\geq4$ and $m\equiv1\mod3$ from Theorem \ref{theorem12} as well. See figure 1.

\begin{figure}[H]
\includegraphics[width=.8\linewidth]{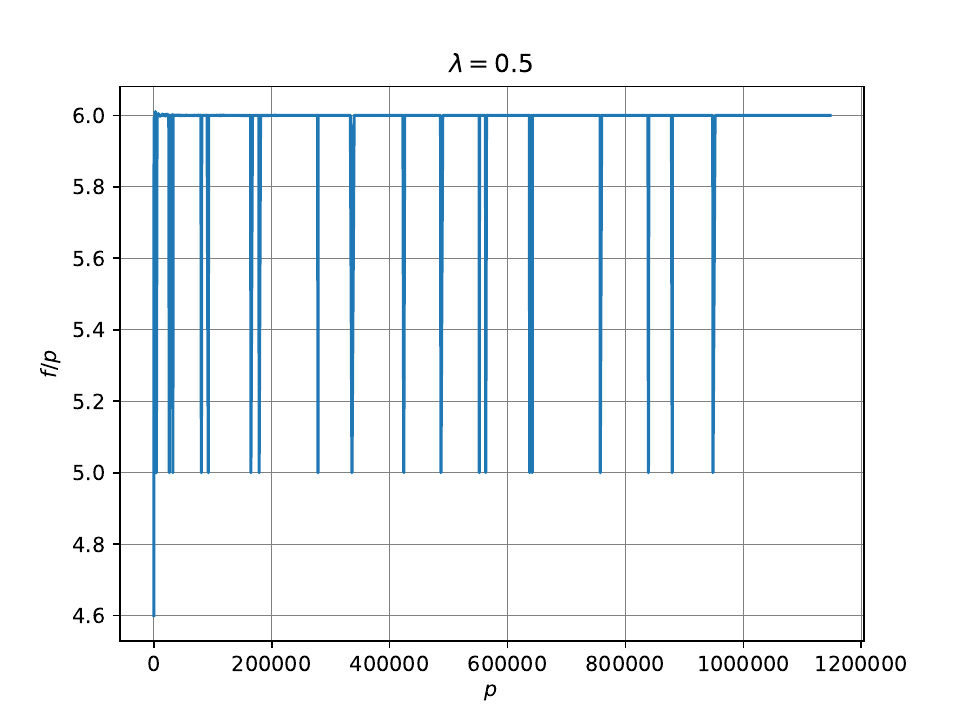}
\caption*{Figure 1: Plot of $f_{1/2}(p)/p$}\end{figure}

By computational evidence (see \cite{repo}), we suspect that $\lim_{p\to\infty}\frac{f_\lambda(p)}p=3$ for $\lambda>1$. See the algorithm in \cite{repo} and figure 2 below.

\begin{figure}[H]
\includegraphics[width=.8\linewidth]{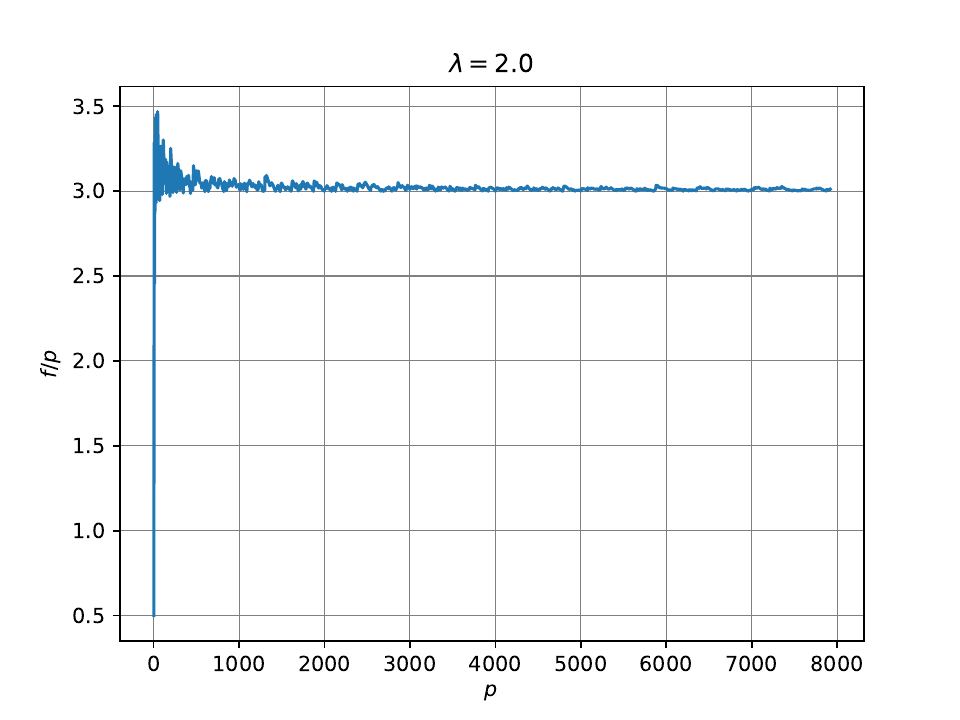}
\caption*{Figure 2: Plot of $f_2(p)/p$}\end{figure}

In particular for the Frobenius number $f_n$ of $S_n$ we would have $\lim_{n\to\infty}\frac{f_n}{p_n}=3$, hence large even numbers would be the sum of two primes (cf. \cite[Proposition~2]{hellus_rechenauer_waldi_19}).

\section{Bounds for the Frobenius number of certain numerical semigroups generated by primes}\label{section_upper_bounds}

In order to verify Theorem~\ref{theorem12}, for integers $m\geq2$ consider the statement

\vspace{.1cm}

\noindent$\mathbf{A(m)}$ For every $\delta>0$ there is an $N(\delta,m)>0$ such that all integers $N\geq N(\delta,m)$ of the same parity as $m$ can be written as a sum of primes
\[\tag{2}N=q_1+\dots+q_m\text{ with the restriction }\left|\frac Nm-q_i\right|<\delta\cdot N\text{ for }i=1,\ldots,m.\]
For short: ``Large $N$ of the same parity as $m$ are sums of $m$ almost equal primes.''

\begin{prop}\label{prop21}

Suppose $A(m)$ holds for some $m\geq2$. Then for each pair $(\varepsilon,\lambda)$ of reals $\varepsilon>0$ and $\lambda>\frac2m$,
\[\tag{3}f_\lambda(p)<(m+1+\varepsilon)p\text{ for large }p\in D_\lambda.\]

\end{prop}

\noindent\textit{Proof:} W.\,l.\,o.\,g. we may assume that $\varepsilon<\lambda-\frac2m$. Let $p\in D_\lambda$. Set $\delta:=\frac\varepsilon{m(m+2+\varepsilon)}$. From $A(m)$ we get for large $p$:

Every integer $N\in[(m+\varepsilon)p,(m+2+\varepsilon)p]$ of the same parity as $m$ can be written as $N=q_1+\dots+q_m$ with primes $q_i$ such that
\[\tag{4}\left|\frac Nm-q_i\right|<\delta\cdot N\text{ for }i=1,\ldots,m.\]
We will see in a moment, that $q_i\in S_\lambda(p)$ for $i=1,\ldots,m$, hence $N\in S_\lambda(p)$: Inequality (4) implies
\[q_i>\frac Nm-\delta\cdot N=\frac{m+2}{m(m+2+\varepsilon)}N\geq\frac{(m+2)(m+\varepsilon)}{m(m+2+\varepsilon)}p>p\]
and
\[q_i<\frac Nm+\delta\cdot N=\frac{m+2+2\varepsilon}{m(m+2+\varepsilon)}N\leq\frac{m+2+2\varepsilon}mp\leq p+\left(\frac2m+\varepsilon\right)p<p+\lambda p,\]
since $m\geq2$ and $\varepsilon<\lambda-\frac2m$. Hence $N\in S_\lambda(p)$.

Considering $N$ and $N+p$ we see, that $\mathbb Z\cap[(m+1+\varepsilon)p,(m+2+\varepsilon)p]$ is contained in $S_\lambda(p)$ and is a set of at least $p$ consecutive integers. Hence $f_\lambda(p)<(m+1+\varepsilon)p$.\hfill$\qedsymbol$

It will be immediate from the following improved version of \cite[Theorem 1.1]{matomaki_etal_17}, that $A(m)$ holds for all integers $m\geq3$.

\begin{prop}

\label{prop22}Let $\theta=\frac{11}{20}+\varepsilon$, $\varepsilon>0$, and $m\geq3$ an integer. Then every sufficiently large integer $N$ of the same parity as $m$ can be written as the sum $N=q_1+\dots+q_m$ of $m$ primes with
\[\left|\frac Nm-q_i\right|\leq N^\theta\text{ for }i=1, \ldots,m.\]

\end{prop}

\noindent\textit{Proof,} due to Kaisa Matom\"aki (private communication \cite{matomaki_2020}).

Given $\varepsilon>0$ and $m\geq3$ as above, let $N$ be a sufficiently large integer with $N-m$ even. By the existence of primes in short intervals, cf. \cite[Theorem 1]{baker_etal_2001}, there is a prime $p\in\left[\frac Nm-\left(\frac Nm\right)^\frac{21}{40},\frac Nm\right]$, hence
\[\tag{5}0\leq\frac Nm-p\leq\left(\frac Nm\right)^{\frac{21}{40}}<N^{\frac{11}{20}+\frac\varepsilon2}<N^{\frac{11}{20}+\varepsilon}.\]
Now use \cite[Theorem 1.1]{matomaki_etal_17} for the odd integer $n:=N-(m-3)p=3\frac Nm+(m-3)\left(\frac Nm-p\right)$ and $\theta=\frac{11}{20}+\frac\varepsilon2$. You get primes $q_1, q_2, q_3$ with $n=q_1+q_2+q_3$ and
\[\tag{6}\left|\frac n3-q_i\right|\leq n^{\frac{11}{20}+\frac\varepsilon2}\leq N^{\frac{11}{20}+\frac\varepsilon2}\text{ for }i=1,2,3.\]
For $N\gg0$ we also have
\[\tag{7}\frac m3\leq N^{\frac\varepsilon2}.\]
Finally by (5), (6) and (7)
\begin{align*}\left|\frac Nm-q_i\right|&=\left|\left(\frac n3-q_i\right)-(m-3)\frac{\frac Nm-p}3\right|\\&\leq N^{\frac{11}{20}+\frac\varepsilon2}+\frac{m-3}3N^{\frac{11}{20}+\frac\varepsilon2}\\&=\frac m3N^{\frac{11}{20}+\frac\varepsilon2}\\&\leq N^{\frac{11}{20}+\varepsilon}\text{ for }i=1, 2, 3.\end{align*}
Hence $n=q_1+q_2+q_3+(m-3)p$ is the sum of $m$ primes of size as desired.\hfill$\qedsymbol$

\begin{coro}

\label{coro23}$A(m)$ is true for all integers $m\geq3$.

\end{coro}

\noindent\textit{Proof.} Let $m\geq3$ and $\delta>0$. Apply \ref{prop22} with $\varepsilon=\frac1{20}$. Then $1>\theta=\frac35>\frac{11}{20}$ and $N^\theta<\delta\cdot N$ for large $N$. Hence $A(m)$ is true by \ref{prop22}.\hfill\qedsymbol

\vspace{.1cm}

Proposition~\ref{prop21} and Corollary \ref{coro23} together imply the following result.

\begin{coro}

\label{coro24}Let $\varepsilon>0$, $m\geq3$ and $\lambda>\frac2m$. Then

\[\pushQED{\qed}f_\lambda(p)<(m+1+\varepsilon)p\text{ for large }p\in D_\lambda.\qedhere\popQED\]

\end{coro}

\noindent\textbf{Remark.} Let $p=p_n$ the $n$th prime in the natural order, $S_n$ the semigroup generated by all primes not less than $p$ and $f_n$ its Frobenius number. $S_{\lambda=1}(p_n)$ is contained in $S_n$, hence $f_n$ is at most $f_{\lambda=1}(p_n)$. Therefore, an application of \ref{coro24} with $m=3$ gives $\limsup_{n\to\infty}\frac{f_n}{p_n}\leq4$. This has been shown by a somewhat different way in our former paper \cite[Remark~2.\,a)]{hellus_rechenauer_waldi_19}.

\vspace{.1cm}

On our way to Theorem~\ref{theorem12}, Proposition \ref{prop25} below will give us lower bounds for $\frac{f_\lambda(p)}p$, $p\in D_\lambda$. Let $p(\lambda):=\max\left(I_\lambda(p)\cap\mathbb P\right)$. From now on let $m\geq2$. Set
\[T(m):=\{t\in\mathbb N\mid1+tm, 3+tm, 1+t(m+2)\in\mathbb P\}.\]

\begin{lemma}\label{lem26}

\begin{enumerate}

\item[a)] Let $p\in D_{\frac2m}$ and $p>m$. If $f_{\frac2m}(p)<(m+2)p-2$, then there is a $t\in T(m)$ such that $p=1+tm$, in particular $p+2$ is a prime as well.

\item[b)] If $T(m)$ is finite, then
\[\frac{f_{\frac2m}(p)}p\geq m+2-\frac2p\text{ for large }p\in D_{\frac2m}.\]

\item[c)] $T(2)=\{1\}$, and $T(m)$ is empty if $m>2$ and $m$ is incongruent to $1$ modulo $3$.

\end{enumerate}

\end{lemma}

\noindent\textit{Proof.} Let $p\in D_\lambda$.

\begin{enumerate}

\item[a)] By definition of $p(\frac2m)$, $mp(\frac2m)\leq (m+2)p$, hence $mp(\frac2m)<(m+2)p$ since $p>m$. For reasons of parity we even have $mp(\frac2m)\leq(m+2)p-2$.

$z:=(m+2)p+2>f_{\frac2m}(p)$, hence $z\in S_{\frac2m}(p)$. Since $mp(\frac2m)<z<(m+3)p$, because of parity $z$ is the sum of exactly $m+2$ atoms from $S_{\frac2m}(p)$, hence $p+2$ must be a prime. Similarly $w:=(m+2)p-2>f_{\frac2m}(p)$ is in $S_{\frac2m}(p)$; hence $w=mp(\frac2m)$ because of its parity and since $mp(\frac2m)\leq w<(m+2)p$. For $t:=\frac{p-1}m$ we have primes $p=1+tm$, $p+2=3+tm$ and $p(\frac2m)=1+t(m+2)$. Further $2t=p(\frac2m)-p$ is an even integer; hence $t\in T(m)$.

\item[b)] is immediate from a).

\item[c)] Let $m\geq2$ and $t>0$ be integers. Elementary calculations modulo $3$ show: If $m$ is incongruent to $1$ modulo $3$, then $3$ divides $(1+tm)(3+tm)(1+t(m+2))$. Hence $1+tm$, $3+tm$ and $1+t(m+2)$ are primes if and only if $m=2$ and $t=1$.\hfill$\qedsymbol$

\end{enumerate}

\vspace{.1cm}

\begin{prop}\label{prop25}

\begin{enumerate}

\item[a)] $\frac{f_\lambda(p)}p\geq3-\frac6p$ for all $\lambda>0$ and $p\in D_\lambda$ (cf. \cite[Proposition~1]{hellus_rechenauer_waldi_19}).

\item[b)] If $m$ is incongruent to $1$ modulo $3$, then
\[\tag{8}\frac{f_{\frac2m}(p)}p\geq m+2-\frac2p\text{ for large }p\in D_{\frac2m}.\]
In case $m\equiv1\mod3$, (8) at least holds for large isolated primes.

\item[c)] Let $m\geq2$ be arbitrary and $0<\lambda<\frac2m$. Then
\[\frac{f_\lambda(p)}p\geq m+2-\frac2p\text{ for }p>\frac2{2-\lambda\cdot m}, p\in D_\lambda.\]

\end{enumerate}

\end{prop}

\noindent\textit{Proof:}

\begin{enumerate}

\item[b)] immediately follows from Lemma~\ref{lem26}.

\item[c)] Elementary calculation shows that, since $p>\frac2{2-\lambda\cdot m}$,
\[m(1+\lambda)p<(m+2)p-2,\]
consequently
\[m\cdot p(\lambda)\leq m(1+\lambda)\cdot p<(m+2)p-2<(m+2)p.\]
Hence $(m+2)p-2$ is a gap of $S_\lambda(p)$, for reasons of parity and magnitude.

\hfill$\qedsymbol$
\end{enumerate}

Now we are ready to restate and prove:

\vspace{.1cm}

\noindent\textbf{Theorem 1.2.} Let $0<\lambda\leq1$, $\lambda\notin\{\frac24,\frac27,\frac2{10},\frac2{13},\ldots\}$, be a real. Let $D_\lambda:=\{p\in\mathbb P\mid[p,p+\lambda p]\text{ contains at least two primes}\}$ and let $f_\lambda(p)$ be the Frobenius number of the numerical semigroup generated by all primes from $[p,p+\lambda p]$. Then 
\[\lim_{p\in D_\lambda\atop p\to\infty}\frac{f_\lambda(p)}p=\left\lfloor2+\frac2\lambda\right\rfloor.\]
For $\lambda=\frac2m$ with an integer $m\geq2$ and $m\equiv1\mod3$, at least
\[\pushQED{\qed}\limsup_{p\to\infty}\frac{f_\lambda(p)}p=\left\lfloor2+\frac2\lambda\right\rfloor.\]

\noindent\textit{Proof.} Immediate from Corollary~\ref{coro24} and Proposition~\ref{prop25}.\hfill$\qedsymbol$

\vspace{.1cm}

Notice, that in case $m\equiv1\mod3$ Dickson's conjecture \cite{dickson_1904} implies, that $T(m)$ is infinite, so the above proof probably will not work. So in this case, Lemma \ref{lem26} does not include formula (8) for large $p$.

\section{The question of Wilf for certain numerical semigroups}
\label{section_Wilf_certain_num_sem_grps}

This section is devoted to the proof of Proposition~\ref{prop11}, restated below for the reader's comfort.

\vspace{.1cm}

Since Wilf's inequality (1) holds by \cite{bruns_etal_2019} and \cite{eliahou_2018} if $p<19$ or $f<3p$, in what follows we may assume that $p\geq19$ and $f>3p$. Proposition \ref{prop11}\,a) will be an easy consequence of the following result.

\begin{theorem}

\label{theorem31}There is a constant $c>0$ such that every numerical semigroup $S$ of multiplicity $p\geq c$ ($p$ not necessarily a prime number)  and containing the primes from $J(p):=[p,p+p^{0.525}]$, satisfies Wilf's inequality
\[\tag{9}\frac g{1+f}\leq1-\frac1e\text{, equivalently }e(1+f-g)\geq1+f.\]

\end{theorem}

\textit{Proof.} Let $\pi(x)$ be the number of primes less than or equal to $x$. For sufficiently large integers $p$, by \cite[p.\,562]{baker_etal_2001} we have
\[\tag{10}e\geq\left|J(p)\cap\mathbb P\right|\geq0.09\cdot\frac{p^{0.525}}{\log p}=:e^*(p)\geq2p^{0.5},\]
the latter since $\lim_{x\to\infty}\frac{x^{0.025}}{\log x}=\infty$. Notice, that $1+f-g$ is the number of elements of $S$ lying below $f$, sometimes called \textit{sporadic} for $S$. Let $m$ be an integer such that $mp<f<(m+1)p$. We have $m\geq3$ since by assumption $f>3p$. Hence the $(m-1)\cdot\left|J(p)\cap\mathbb P\right|$ many elements
\[s=ip+q, 0\leq i\leq m-2\text{ and }q\in J(p)\text{ a prime},\]
are sporadic for $S$, since $s\leq(m-2)p+q<mp<f$. Finally we get by (10) and since $m\geq3$
\[\pushQED{\qed}e(1+f-g)\geq e^*(p)\cdot(m-1)\cdot e^*(p)\geq(m-1)4p>(m+1)p>f.\qedhere\popQED\]

\vspace{.1cm}

\begin{lemma}

\label{lemma333}Let $S(p)$ be the semigroup generated by the primes from $I(p)=[p,2p]$ and $f(p)$ its Frobenius number. Then
\[\tag{11}f(p)<2(\pi(2p)-\pi(p))^2\text{, if }n=\pi(p)>674.\]

\end{lemma}

\noindent\textit{Proof.} Fundamental for this are the approximate formulas for the functions $p_n$ and $\pi(x)$ from the papers \cite{rosser_schoenfeld_1962} and \cite{rosser_schoenfeld_1975} by Rosser and Schoenfeld. According to \cite{rosser_schoenfeld_1975} we have
\[\tag{12}\pi(2x)<2\pi(x)\text{ for }x\geq11.\]
In \cite{rosser_schoenfeld_1962} it is shown, that
\[\tag{13}p_n<n(\log n+\log\log n)\text{ for }n>5,\]
\[\tag{14}\pi(x)<\frac x{\log x-\frac32}\text{ if }\log x>\frac32\text{ and}\]
\[\tag{15}\pi(x)>\frac x{\log x-\frac12}\text{ for }x\geq67.\]
Since for the embedding dimension $e(p_n)$ of $S(p_n)$ we have
\[e(p_n)=\pi(2p_n)-n+1<n+1<p_n\text{ for }n>674\text{ by (12),}\]
the approximation of the Frobenius number by \cite{selmer_1977} page 2, last line can be applied to $S(p)$ if $\pi(p)=n>674$. We get:
\begin{align*}\tag{16}f(p)&\leq2 p_{\pi(2p)}\cdot\lfloor p/(\pi(2p)-n +1)\rfloor-p\\&<2pp_{\pi(2p)}/(\pi(2p)-n+1)\\&<2p_n\cdot p_{2n}/(\pi(2p)-n).\end{align*}
From (14) and (15) we get
\[\tag{17}\pi(2x)>2\cdot\frac x{\log(2x)-\frac12}>2\frac{\log x-\frac32}{\log(2x)-\frac12}\cdot\pi(x)=:l(x)\cdot\pi(x)\text{, }x\geq67.\]
It is easily seen that the function $l(x)$, $x\geq67$, is strictly increasing. Together with (12) and (17) we get at the places $x=p_n$, $n\geq675$, i.\,e. $p_n\geq5039$
\[\tag{18}2n>\pi(2p_n)>l(p_n)\cdot\pi(p_n)\geq l(5039)\cdot n,\]
where $l(5039)$ is approximately $1.61158$.

The function
\[l_2(x):=\frac{2\cdot(\log x+\log\log x)(\log(2x)+\log\log(2x))}x\text{, }x\geq675\]
is strictly decreasing. As one can check,
\[\tag{19}l_2(675)<(l(5039)-1)^3.\]
Applying (13) to the right hand side of formula (16), we get for $n\geq675$
{\belowdisplayskip=-12pt\begin{align*}f(p_n)&<2p_np_{2n}/(\pi(2p_n)-n)\\&\buildrel\text{(13)}\over<2\cdot l_2(n)\frac{n^3}{\pi(2p_n)-n}\\&\leq2\cdot l_2(675)\cdot\frac{n^3}{\pi(2p_n)-n}\\&\buildrel\text{(19)}\over<2\cdot(l(5039)-1)^3\cdot\frac{n^3}{\pi(2p_n)-n}\\&\leq2\cdot\frac{((l(p_n)-1)n)^3}{\pi(2p_n)-n}\\&\buildrel\text{(17)}\over<2\cdot(\pi(2p_n)-n)^2.\end{align*}\hfill$\qedsymbol$}

\vspace{.2cm}

For $\lambda>0$, let $D_\lambda:=\{p\in\mathbb P\mid[p,p+\lambda p]\text{ contains at least two primes}\}$.

\vspace{.1cm}

\noindent\textbf{Proposition 1.1.} Let $\lambda>0$, $p\in D_\lambda$ and $S$ any numerical semigroup of multiplicity $p$ containing $[p,p+\lambda p]\cap\mathbb P$.

\begin{enumerate}

\item[a)] There is an integer $C(\lambda)>0$ such that for $p>C(\lambda)$, the semigroups $S$ from above satisfy Wilf's inequality (1).

\item[b)] In case $\lambda=1$ formula (1) holds for all p.
In particular (1) is true for $S=S_n:=$numerical semigroup generated by all primes not less than $p_n$.

\end{enumerate}

\noindent\textit{Proof.} a) Let $\lambda>0$, $p\in D_\lambda$ and $S$ as in the statement. Let $c$ be the constant from Theorem \ref{theorem31}, choose $C(\lambda)\geq c$ such that $C(\lambda)^{0.525}<\lambda\cdot C(\lambda)$. Then for every prime $p\geq C(\lambda)$ we have $p^{0.525}<\lambda p$ as well. Hence $p\geq c$ and $J(p)\cap\mathbb P\subseteq[p,p+\lambda p]\cap\mathbb P\subseteq S$, as requested in Theorem~\ref{theorem31}, and $S$ satisfies (9).

b) Let $S$ be as in \ref{prop11}\,b), $S(p)$ the semigroup generated by the primes from $I(p):=[p,2p]$ and $f(p)$ its Frobenius number. Since $3p<f$, the primes from $I(p)\subseteq S$ are atoms as well as sporadic elements for $S$. The latter also holds for the even numbers $p+q$, $q$ a prime from $I(p)$ as well as for $3p$. Hence $1+f-g\geq2(\pi(2p)-\pi(p)+1)+1$, and all together
\[\tag{20}e(1+f-g)\geq2(\pi(2p)-\pi(p)+1)^2+\pi(2p)-\pi(p)+1.\]
Since $f\leq f(p)$, (20) together with (11) from Lemma \ref{lemma333} imply (9) for $n>674$. Therefore, it remains only to prove Proposition~\ref{prop11}\,b) in case $7<n<675$.

According to the last column of table \texttt{wilf\_for\_p\_to\_2p.pdf} from \cite{repo} inequality (9) holds for $S(p)$, if $8\leq\pi(p)\leq675$.

Hence we may assume that $S$ is different from $S(p)$, $p=p_n$. Then $e\geq\pi(2p_n)-n+2$, and (20) can be improved to
\[e\cdot(1+f-g)\geq(2\cdot(\pi(2p_n)-n+1)+1)(\pi(2p_n)-n+2).\]
The second last column of table \texttt{wilf\_for\_p\_to\_2p.pdf} from \cite{repo} mentioned above shows, that $f(p)$ (and, a fortiori, $f$) is less than the right hand side of this inequality, if $10\leq\pi(p)<675$. The remaining cases are $p=19$ and $p=23$.

For $p=23$, by assumption we have $f>69 = 3\cdot p$, and $S(23)$ contains $17$ elements less than $70$, which then are sporadic for $S$. Since $S$ is different from $S(23)$, we have $e\geq e(23)+1=7$; finally $e(1+f-g)\geq7\cdot17>102=f(23)\geq f$.

Analogously for $p=19$ we have $f>57 = 3\cdot p$ and $e\geq6$. Further $58=29+29$ is in $S$, hence $f\geq59$. Since $60$, $61$ and $62$ are also in $S$, either $f=59$ or $f\geq63$.

\begin{enumerate}

\item[a)] Case $f\geq63$: $S(19)$ contains $19$ elements $<63$. As above $e(1+f-g)\geq6\cdot19>101=f(19)\geq f$.

\item[b)] Case $f=59$: $S(19)$ contains $16$ elements $<59$. It follows $e(1+f-g)\geq6\cdot16>59=f$.\hfill$\qedsymbol$

\end{enumerate}

The fraction $d=\frac{1+f-g}{1+f}$ describes the density of the sporadic elements of $S$  in $[0,f]\cap \mathbb Z$. In terms of this density, Wilf's conjecture says that $d$ is at least $\frac1e$.

We will see in a moment that for the semigroups $S(p)$ generated by the primes from $[p,2p]$ this bound becomes extremely weak, as $p$ goes to infinity.

As above let $2\geq\lambda>0$ be a real parameter, $S_\lambda(p)$ the semigroup generated by the primes from $I_\lambda(p)$ and $e_\lambda(p)$ its embedding dimension. Here $e_\lambda(p)\sim\lambda\cdot\pi(p)$, hence $\frac1{e_\lambda(p)}$ is a null sequence. In contrast, the results of \cite{coppola_laporta_95} and \cite{matomaki_etal_17} imply the following result.

\begin{prop}

If $S=S(p)$ is the semigroup generated by the primes from $[p,2p]$ then
\[d\sim\frac38.\]

\end{prop}

\noindent\textit{Proof} Let $\frac58<t<1$. By \cite[Corollary]{coppola_laporta_95}, for all $2N\in[2p,4p]$ but $O\left(\frac{2p}{\log(2p)}\right)$ exceptions, we have $2N=q_1+q_2$ with
\[N-N^t\leq q_i\leq N+N^t, q_i\text{ a prime for }i=1,2.\]
If even
\[p+(2p)^t\leq N\leq2p-(2p)^t\text{ and  }2N\text{ is not an exception,}\]
then it easily follows that $p<q_i<2p$, hence $2N\in S(p)$. Since $\frac{p^t}p$ and $\frac1{\log(2p)}$ are null sequences, this shows that for large $p$, almost all even elements from $[2p,4p]$ are in $S(p)$.

By similar arguments we see from \cite[Theorem~1.1]{matomaki_etal_17}, that for large primes $p$, every odd integer $N$ with
\[p+(6p)^t\leq\frac N3\leq2p-(6p)^t\]
is contained in $S(p)$. Hence for large $p$, almost all odd elements from $[3p,6p]$ are in $S(p)$. Further by Theorem \ref{theorem12}, $f(p)\sim4p$.\hfill$\qedsymbol$

\section{Binary Goldbach for large numbers: Sufficient conditions}

The \textit{Binary Goldbach conjecture} for large numbers, which seems to be open, states that each large enough even integer can be written as a sum of two primes.

In this section we present some consequences which would follow if ``Binary Goldbach for large numbers'' should be false. We cannot disprove any of these consequences and we do not believe that our results mean some practical progress on a way to prove Binary Goldbach for large numbers.

Recall that by $S_n$ we denote the semigroup generated by all primes not less than $p_n$, and by $f_n$ we denote the Frobenius number of $S_n$.

Obviously $3p_n\notin S_{n+1}$, hence $3p_n\leq f_{n+1}$ for all $n\geq2$.

On the other hand, the table \texttt{full\_numerical\_semigroups\_created\_by}
\texttt{primes\_greater\_nth.pdf} from \cite{repo} shows that, for $2\leq n\leq10\, 000$,

\begin{enumerate}

\item[(i)] $f_{n+1}<3p_n+2n$.

\item[(ii)] $f_{n+1}$ is odd, with the exception $f_4=16$.

\item[(iii)] $f_{n+1}=3p_n+2n-2$ for $n=4,6,7,9$ or $15$.

\end{enumerate}

Analogously, improving \cite[Lemma~3]{hellus_rechenauer_waldi_19} we get the following result.

\begin{lemma}

\label{lemma41}For large $n$, each odd integer $N\geq 3p_n+2n$ is contained in $S_{n+1}$. In particular,
\[f_{n+1}<3p_n+2n\text{ if }f_{n+1}\text{ is odd, and}\]
\[f_{n+1}<3p_n+p_{n+1}+2n\text{ if }f_{n+1}\text{ is even.}\]

\end{lemma}

\noindent\textit{Proof: }  By \cite[Theorem~1.1]{matomaki_etal_17} each odd number $N\geq U_n:=3p_n+2n\gg0$ is the sum
\[N=q_1+q_2+q_3\text{ of primes }q_i\text{ with }\left|\frac N3-q_i\right|\leq N^{\frac35}\text{, }i=1,2,3.\]
Since $F(x):=\frac x3-x^{\frac35}$ is increasing for $x\geq5$, and $p_n\sim n\log n$,

\begin{align*}q_i&\geq\frac N3-N^{\frac35}\\&\geq\frac{U_n}3-U_n^{\frac35}\\&>p_n+\frac{2n}3-(6n\log n +2n)^{\frac35}\\&>p_n\end{align*}

for large $n$. Hence $N=q_1+q_2+q_3$ is contained in $S_{n+1}$.\hfill$\qedsymbol$

\vspace{.2cm}

This implies for each $1\geq \varepsilon>0$: Since $\frac{2n}{p_n}$ is a null sequence and $\lim_{n\to\infty}\frac{p_{n+1}}{p_n}=1$, there is an $n(\varepsilon)>0$ such that for $n\geq n(\varepsilon)$

\[\tag{21}\frac{f_{n+1}}{p_n}<3+\varepsilon\text{ if }f_{n+1}\text{ is odd, }\frac{f_{n+1}}{p_{n+1}}<4+\varepsilon\text{ if }f_{n+1}\text{ is even and }\]
\[\tag{22}\frac{p_{n+1}}{p_n}<1+\frac\varepsilon5.\]
Hence for all $n\geq n(\varepsilon)$
\begin{align*}\tag{23}\frac{f_{n+1}}{p_n}-\frac{f_{n+2}}{p_{n+1}}&\leq\frac{f_{n+1}}{p_{n+1}}\cdot\frac{p_{n+1}}{p_n}-\frac{f_{n+1}}{p_{n+1}}\\&=\frac{f_{n+1}}{p_{n+1}}\left(\frac{p_{n+1}}{p_n}-1\right)\\&<5\cdot\frac\varepsilon5\\&=\varepsilon\end{align*}
by (21) and (22).

\noindent\textbf{Assumption} The Binary Goldbach conjecture for large numbers is false.

\vspace{.1cm}

\noindent\textbf{Conclusions for the sequence }$\mathbf{(f_n)}$

\begin{itemize}

\item [a)] For infinitely many $n>0$, $\frac{f_{n+1}}{p_n}\geq4$. See \cite[Lemma~1]{hellus_rechenauer_waldi_19}.

\item [b)] For each integer $k>0$ there is an integer $n\geq n\left(\frac1k\right)$ such that the $k$ Frobenius numbers $f_{n+1}, \ldots, f_{n+k}$ are even.

\noindent\textit{Proof:} By a) we can find an $n\geq n\left(\frac1k\right)$ such that $\frac{f_{n+1}}{p_n}\geq4$. By (23), for $1\leq m\leq k$ we have

\begin{align*}\frac{f_{n+m}}{p_{n+m-1}}&=\frac{f_{n+1}}{p_n}-\left(\frac{f_{n+1}}{p_n}-\frac{f_{n+2}}{p_{n+1}}\right)-\dots-\left(\frac{f_{n+m-1}}{p_{n+m-2}}-\frac{f_{n+m}}{p_{n+m-1}}\right)\\&\geq4-(m-1)\cdot\frac1k\\&\geq3+\frac1k,\end{align*}

hence the integers $f_{n+1}, \ldots, f_{n+k}$ are even by (21).\hfill$\qedsymbol$

\item [c)] Either $f_{n+1}$ is even for almost all $n>0$, or $[3,4]$ is contained in the closure of $\left\{\left.\frac{f_{n+1}}{p_n}\right|n>0\right\}$.

\noindent\textit{Proof:} Under the \textbf{additional assumption}, that $f_{n+1}$ is odd for infinitely many $n$ we have to show:

Let $1\geq\varepsilon>0$ and $x\in[3,4]$ be arbitrary. Then
\[\left|x-\frac{f_{m+1}}{p_m}\right|<\varepsilon\text{ for some integer }m>0.\]
\noindent\textit{Proof:} By a) and since $f_{n+1}$ is odd infinitely often, there are integers $n\geq n(\varepsilon)$ and $k>0$ such that
\[\frac{f_{n+1}}{p_n}\geq4\text{ and }\frac{f_{n+k+1}}{p_{n+k}}<3+\varepsilon.\]
In case $x\leq\frac{f_{n+k+1}}{p_{n+k}}$ we take $m=n+k$. Otherwise
\[\frac{f_{n+k+1}}{p_{n+k}}<x\leq4\leq\frac{f_{n+1}}{p_n},\]
and we can find an $m$ with $n+1\leq m\leq n+k$ such that $\frac{f_{m+1}}{p_m}<x\leq\frac{f_m}{p_{m-1}}$.

Hence by (23)
\[\pushQED{\qed}0\leq\frac{f_m}{p_{m-1}}-x<\frac{f_m}{p_{m-1}}-\frac{f_{m+1}}{p_m}<\varepsilon.\qedhere\popQED\]

\end{itemize}

\noindent\textbf{Note} on the even $f_n$:

\begin{enumerate}

\item We apply \cite[Corollary]{coppola_laporta_95} to any $0<\varepsilon<\frac38$ like e.\,g. $\varepsilon=\frac18$: $U=m^{\frac58+\varepsilon}$ is, because of $\varepsilon<\frac38$, smaller than $\frac m{10}$ for $m$ large.

Hence apart from at most $O\left(\frac N{(\log N)^A}\right)$ exceptions ($A$ arbitrary) each even integer
\[2m\in[N,2N]\]
is a sum
\[2m=q_1+q_2\]
of two primes numbers
\[\frac9{10}m\leq q_1,q_2\leq\frac{11}{10}m.\]

\item Let $n$ be large enough and suppose $f_n$ is even. By \cite[Proposition~1 and Lemma~3]{hellus_rechenauer_waldi_19}, $N:=3p_n-6\leq f_n\leq 2N$; in particular, $p_n\leq\frac9{10}\cdot\frac{f_n}2$. Hence the gap $f_n$ of $S_n$ is always an exception in the sense of (1).

\end{enumerate}

\textbf{Note on coding} Our numerical experiments may be reproduced by using the corresponding codes from the repository \texttt{Bleiglanz/On\_The\_FrobeniusNum\-ber} in \cite{repo}.

\section*{Acknowledgement}

The authors are grateful to the referee for comments and suggestions and  would also like to acknowledge the creators of the following used software: GAP (\cite{GAP}), NumericalSgps (\cite{DGM19}), Python (\cite{python}), Matplotlib (\cite{matplotlib}), NumPy (\cite{numpy}), Rust (\cite{rust}).

\end{document}